\documentclass[12pt,eqno,]{article}
\usepackage{amsmath,amsthm,amssymb}

\makeatletter

\newcommand{\Rmnum}[1]{\expandafter\@slowromancap\romannumeral #1@}
\makeatother

\newcommand{\arcsinh}{\rm {arcsinh}}
\title{\normalsize\bf  REFINEMENTS OF THE INEQUALITIES BETWEEN NEUMAN-S\'{A}NDOR, ARITHMETIC, CONTRA-HARMONIC AND QUADRATIC MEANS}
\author{\small  YU-MING CHU AND MIAO-KUN WANG}
\date{}

\begin{document}
\maketitle

\renewcommand{\thefootnote}{\fnsymbol{footnote}}

\footnotetext{\hspace*{-5mm}
\begin{tabular}{@{}r@{}p{14.0cm}@{}}
&\qquad Mathematics Subject Classification (2010): 26E20.\\
&\qquad Keywords and phrases: Neuman-S\'{a}nder mean, arithmetic mean, contra-harmonic mean, quadratic mean.\\
&\qquad This research was supported by the Natural Science
Foundation of China under Grants 11071069 and 11171307, and Innovation
Team Foundation of the Department of Education of Zhejiang Province
under Grant T200924.
\end{tabular}}

\vspace{-1cm}

\def\abstractname{}
\begin{abstract}
\noindent {\emph{Abstract}.} In this paper, we prove that the
inequalities $\alpha [1/3 Q(a,b)+2/3
A(a,b)]+(1-\alpha)Q^{1/3}(a,b)A^{2/3}(a,b)<M(a,b) <\beta [1/3
Q(a,b)+2/3 A(a,b)]+(1-\beta)Q^{1/3}(a,b)A^{2/3}(a,b)$ and $\lambda
[1/6 C(a,b)+5/6
A(a,b)]+(1-\lambda)C^{1/6}(a,b)A^{5/6}(a,b)<M(a,b)<\mu [1/6
C(a,b)+5/6 A(a,b)]++(1-\mu)C^{1/6}(a,b)A^{5/6}(a,b)$ hold for all
$a,b>0$ with $a\neq b$ if and only if $\alpha\leq
(3-3\sqrt[6]{2}\log(1+\sqrt{2}))/[(2+\sqrt{2}-3\sqrt[6]{2})\log(1+\sqrt{2})]=0.777\cdots$,
$\beta\geq 4/5$, $\lambda\leq
(6-6\sqrt[6]{2}\log(1+\sqrt{2}))/(7-6\sqrt[6]{2}\log(1+\sqrt{2}))=0.274\cdots$,
and $\mu\geq 8/25$. Here, $M(a,b)$, $A(a,b)$, $C(a,b)$, and $Q(a,b)$
denote the Neuman-S\'{a}ndor, arithmetic, contra-harmonic, and
quadratic means of $a$ and $b$, respectively.
\end{abstract}

\bigskip
\centerline{\bf 1. Introduction} \setcounter{section}{1}
\setcounter{equation}{0}
\bigskip

For $a,b>0$ with $a\neq b$ the Neuman-S\'{a}ndor mean $M(a,b)$ [1] is defined by
\begin{equation}
M(a,b)=\frac{a-b}{2{{\arcsinh}}\left[(a-b)/(a+b)\right]}.
\end{equation}

Recently, the Neuman-S\'{a}ndor mean has been the subject intensive research. In particular, many remarkable inequalities for the  Neuman-S\'{a}ndor mean $M(a,b)$ can be found in the literature [1, 2].

Let $A(a,b)=(a+b)/2$, $G(a,b)=\sqrt{ab}$, $L(a,b)=(b-a)/(\log{b}-\log{a})$, $C(a,b)=(a^2+b^2)/(a+b)$, $Q(a,b)=\sqrt{(a^2+b^2)/2}$, $P(a,b)=(a-b)/(4\arctan{\sqrt{a/b}}-\pi)$ and $T(a,b)=(a-b)/\left[2\arcsin((a-b)/(a+b))\right]$ be the arithmetic, geometric, logarithmic, contra-harmonic, quadratic, first Seiffert and second Seiffert means of $a$ and $b$, respectively. Then it is well known that the inequalities
\begin{equation*}
G(a,b)<L(a,b)<P(a,b)<A(a,b)<T(a,b)<Q(a,b)<C(a,b)
\end{equation*}
hold for all $a,b>0$ with $a\neq b$.

Neuman and S\'{a}ndor [1, 2] established that
\begin{equation*}
A(a,b)<M(a,b)<T(a,b)
\end{equation*}
\begin{equation*}
P(a,b)M(a,b)<A^{2}(a,b)
\end{equation*}
\begin{equation*}
A(a,b)T(a,b)<M^{2}(a,b)<(A^{2}(a,b)+T^{2}(a,b))/2
\end{equation*}
for all $a,b>0$ with $a\neq b$.

Let $0<a,b<1/2$ with $a\neq b$, $a'=1-a$ and $b'=1-b$. Then the following KyFan inequalities
\begin{equation*}
\frac{G(a,b)}{G(a',b')}<\frac{L(a,b)}{L(a',b')}<\frac{P(a,b)}{P(a',b')}<\frac{A(a,b)}{A(a',b')}<\frac{M(a,b)}{M(a',b')}<\frac{T(a,b)}{T(a',b')}
\end{equation*}
were presented in [1].

Very recently, Li et al. [3] proved that
$L_{p_{0}}(a,b)<M(a,b)<L_{2}(a,b)$ for all $a,b>0$ with $a\neq b$,
where $L_{p}(a,b)=[(b^{p+1}-a^{p+1})/((p+1)(b-a))]^{1/p}(p\neq -1,
0)$, $L_{0}(a,b)=1/e(b^{b}/a^{a})^{1/(b-a)}$ and
$L_{-1}(a,b)=(b-a)/(\log{b}-\log{a})$ is the $p$-th generalized
logarithmic mean of $a$ and $b$, and $p_{0}=1.843\cdots$ is the
unique solution of the equation $(p+1)^{1/p}=2\log(1+\sqrt{2})$.
And, in [4] the author found that
\begin{equation}
Q^{1/3}(a,b)A^{2/3}(a,b)<M(a,b)<\frac{1}{3}Q(a,b)+\frac{2}{3}A(a,b)
\end{equation}
and
\begin{equation}
C^{1/6}(a,b)A^{5/6}(a,b)<M(a,b)<\frac{1}{6}C(a,b)+\frac{5}{6}A(a,b)
\end{equation}
for all $a,b>0$ with $a\neq b$.

The aim of this paper is to improve and refine inequalities (1.2) and (1.3). Our main results are the following Theorems 1.1 and 1.2.

\medskip
{\bf THEOREM 1.1.} The double inequality
\begin{align}
\alpha [1/3 Q(a,b)+2/3 A(a,b)]+(1-\alpha)Q^{1/3}(a,b)A^{2/3}(a,b)<M(a,b)\nonumber\\
<\beta [1/3 Q(a,b)+2/3 A(a,b)]+(1-\beta)Q^{1/3}(a,b)A^{2/3}(a,b)
\end{align}
holds for all $a,b>0$ with $a\neq b$ if and only if $\alpha\leq (3-3\sqrt[6]{2}\log(1+\sqrt{2}))/[(2+\sqrt{2}-3\sqrt[6]{2})\log(1+\sqrt{2})]=0.777\cdots$ and $\beta\geq 4/5$.

\medskip
{\bf THEOREM 1.2.} The double inequality
\begin{align}
\lambda [1/6 C(a,b)+5/6 A(a,b)]+(1-\lambda)C^{1/6}(a,b)A^{5/6}(a,b)<M(a,b)\nonumber\\
<\mu [1/6 C(a,b)+5/6 A(a,b)]+(1-\mu)C^{1/6}(a,b)A^{5/6}(a,b)
\end{align}
holds for all $a,b>0$ with $a\neq b$ if and only if $\lambda\leq (6-6\sqrt[6]{2}\log(1+\sqrt{2}))/(7-6\sqrt[6]{2}\log(1+\sqrt{2}))=0.274\cdots$ and $\mu\geq 8/25$.

\medskip
\bigskip
\centerline{\bf 2. Lemmas}
\setcounter{section}{2}\setcounter{equation}{0}
\bigskip

In order to establish our main
results we need two Lemmas, which we present in
this section.

\medskip
{\bf LEMMA 2.1.} Let $p\in(0,1)$,
\begin{equation}
f_{p}(t)={\arcsinh}(\sqrt{t^6-1})-\frac{3\sqrt{t^6-1}}{pt^3+3(1-p)t+2p},
\end{equation}
and $\alpha_{0}=(3-3\sqrt[6]{2}\log(1+\sqrt{2}))/[(2+\sqrt{2}-3\sqrt[6]{2})\log(1+\sqrt{2})]=0.777\cdots$.
Then $f_{4/5}(t)>0$ and $f_{\alpha_{0}}(t)<0$ for all
$t\in(1,\sqrt[6]{2})$.

\medskip
{\em Proof.} From (2.1) one has
\begin{equation}
f_{p}(1)=0,
\end{equation}
\begin{equation}
f_{p}(\sqrt[6]{2})=\ln(1+\sqrt{2})-\frac{3}{p(\sqrt{2}+2)+3\sqrt[6]{2}(1-p)},
\end{equation}
\begin{align}
{f_{p}}'(t)=\frac{3(t-1)^2g_{p}(t)}{[pt^3+3(1-p)t+2p]^2\sqrt{t^6-1}},
\end{align}
where
\begin{align}
g_{p}(t)=&p^2t^6+2p^2t^5+3(-p^2+4p-2)t^4+2(-2p^2+9p-6)t^3\nonumber\\
&+(4p^2+6p-9)t^2-6(1-p)t-3(1-p).
\end{align}

We dived the proof into two cases.

{\bf Case 1} $p=4/5$. Then (2.5) becomes
\begin{equation}
g_{p}(t)=g_{4/5}(t)=\frac{t-1}{25}\left(16t^5+48t^4+90t^3+86t^2+45t+15\right)>0
\end{equation}
for $t\in(1,\sqrt[6]{2})$.

Therefore, $f_{4/5}(t)>0$ for all $t\in(1,2^{1/6})$ follows from
(2.2), (2.4) and (2.6).

{\bf Case 2} $p=\alpha_{0}$. Then (2.3) and (2.5) lead to
\begin{equation}
f_{\alpha_{0}}(\sqrt[6]{2})=0
\end{equation}
and
\begin{align*}
g_{p}(t)=g_{\alpha_{0}}(t)=&\alpha_{0}^2t^6+2\alpha_{0}^2t^5+3(-\alpha_{0}^2+4\alpha_{0}-2)t^4-2(2\alpha_{0}^2-9\alpha_{0}+6)t^3\\
&-(-4\alpha_{0}^2-6\alpha_{0}+9)t^2-6(1-\alpha_{0})t-3(1-\alpha_{0}).
\end{align*}
Note that
\begin{equation}
g_{\alpha_{0}}(1)=9(5\alpha_{0}-4)<0,\quad g_{\alpha_{0}}(\sqrt[6]{2})=0.569\cdots>0,
\end{equation}
\begin{align}
{g_{\alpha_{0}}}'(t)=&6\alpha_{0}^2t^5+10\alpha_{0}^2t^4+12(-\alpha_{0}^2+4\alpha_{0}-2)t^3-6(2\alpha_{0}^2-9\alpha_{0}+6)t^2\nonumber\\
&-2(-4\alpha_{0}^2-6\alpha_{0}+9)t-6(1-\alpha_{0})\nonumber\\
>&6\alpha_{0}^2t^2+10\alpha_{0}^2t^2+12(-\alpha_{0}^2+4\alpha_{0}-2)t^2-6(2\alpha_{0}^2-9\alpha_{0}+6)t^2\nonumber\\
&-2(-4\alpha_{0}^2-6\alpha_{0}+9)t^2-6(1-\alpha_{0})\nonumber\\
=&6(19\alpha_{0}-13)t^2-6(1-\alpha_{0})>12(10\alpha_{0}-7)>0
\end{align}
for $t\in(1,\sqrt[6]{2})$.

The inequality (2.9) implies that
$g_{\alpha_{0}}(t)$ is strictly increasing
in $(1,\sqrt[6]{2})$. Then from (2.4) and (2.8) we clearly see that there exists
$t_{0}\in(1,\sqrt[6]{2})$ such that $f_{\alpha_{0}}(t)$ is strictly decreasing in $[1,t_{0}]$ and strictly increasing in $[t_{0},\sqrt[6]{2}]$.

Therefore, $f_{\alpha_{0}}(t)<0$ for $t\in(1,\sqrt[6]{2})$ follows from (2.2) and (2.7) together with the piecewise monotonicity of $f_{\alpha_{0}}(t)$. $\Box$

\medskip
{\bf LEMMA 2.2.} Let $p\in(0,1)$,
\begin{equation}
F_{p}(t)={\arcsinh}(\sqrt{t^6-1})-\frac{6\sqrt{t^6-1}}{pt^6+6(1-p)t+5p},
\end{equation}
and $\lambda_{0}=6(1-\sqrt[6]{2}\log(1+\sqrt{2}))/(7-6\sqrt[6]{2}\log(1+\sqrt{2}))=0.274\cdots$.
Then $F_{8/25}(t)>0$ and $F_{\lambda_{0}}(t)<0$ for all
$t\in(1,\sqrt[6]{2})$.

\medskip
{\bf Proof.} From (2.10) we have
\begin{equation}
F_{p}(1)=0,
\end{equation}
\begin{equation}
F_{p}(\sqrt[6]{2})=\ln(1+\sqrt{2})-\frac{6}{7p+6\sqrt[6]{2}(1-p)},
\end{equation}
\begin{align}
{F_{p}}'(t)=\frac{3(t-1)^2G_{p}(t)}{[pt^6+6(1-p)t+5p]^2\sqrt{t^6-1}},
\end{align}
where
\begin{align}
G_{p}(t)=&p^2t^{12}+2p^2t^{11}+3p^2t^{10}+2p(3+2p)t^9+p(12+5p)t^8+6p(5-p)t^7\nonumber\\
&+p(48-7p)t^6+2p(33-4p)t^5+3(-3p^2+36p-8)t^4\nonumber\\
&+2(-5p^2+54p-24)t^3+(25p^2+36p-36)t^2\nonumber\\
&-24(1-p)t-12(1-p).
\end{align}

We dived the proof into two cases.

{\bf Case 1} $p=8/25$. Then (2.14) becomes
\begin{align}
&G_{p}(t)=G_{8/25}(t)=\frac{4(t-1)}{625}(16t^{11}+48t^{10}+96t^{9}+460t^{8}+1140t^7\nonumber\\
&+2544t^6+4832t^5+8004t^4+9510t^3+7250t^2+3825t+1275)>0
\end{align}
for $t\in(1,\sqrt[6]{2})$.

Therefore, $F_{8/25}(t)>0$ for all $t\in(1,\sqrt[6]{2})$ follows from
(2.11), (2.13) and (2.15).

{\bf Case 2} $p=\lambda_{0}$. Then (2.12) and (2.14) lead to
\begin{equation}
F_{\lambda_{0}}(\sqrt[6]{2})=0
\end{equation}
and
\begin{align*}
G_{p}(t)=G_{\lambda_{0}}(t)=&\lambda_{0}^2t^{12}+2\lambda_{0}^2t^{11}+3\lambda_{0}^2t^{10}+2\lambda_{0}(3+2\lambda_{0})t^9+\lambda_{0}(12+5\lambda_{0})t^8\nonumber\\
&+6\lambda_{0}(5-\lambda_{0})t^7+\lambda_{0}(48-7\lambda_{0})t^6+2\lambda_{0}(33-4\lambda_{0})t^5\nonumber\\
&+3(-3\lambda_{0}^2+36\lambda_{0}-8)t^4+2(-5\lambda_{0}^2+54\lambda_{0}-24)t^3\nonumber\\
&+(25\lambda_{0}^2+36\lambda_{0}-36)t^2-24(1-\lambda_{0})t-12(1-\lambda_{0}).
\end{align*}
Note that
\begin{equation}
G_{\lambda_{0}}(1)=18(25\lambda_{0}-8)<0,\quad G_{\lambda_{0}}(\sqrt[6]{2})=12.313\cdots>0,
\end{equation}
\begin{align}
{G_{\lambda_{0}}}'(t)=&12\lambda_{0}^2t^{11}+22\lambda_{0}^2t^{10}+30\lambda_{0}^2t^{9}+18\lambda_{0}(3+2\lambda_{0})t^8+8\lambda_{0}(12+5\lambda_{0})t^7\nonumber\\
&+42\lambda_{0}(5-\lambda_{0})t^6+6\lambda_{0}(48-7\lambda_{0})t^5+10\lambda_{0}(33-4\lambda_{0})t^4\nonumber\\
&+12(-3\lambda_{0}^2+36\lambda_{0}-8)t^3-6(5\lambda_{0}^2-54\lambda_{0}+24)t^2\nonumber\\
&-2(-25\lambda_{0}^2-36\lambda_{0}+36)t-24(1-\lambda_{0})\nonumber\\
>&12\lambda_{0}^2t^{2}+22\lambda_{0}^2t^{2}+30\lambda_{0}^2t^{2}+18\lambda_{0}(3+2\lambda_{0})t^2+8\lambda_{0}(12+5\lambda_{0})t^2\nonumber\\
&+42\lambda_{0}(5-\lambda_{0})t^2+6\lambda_{0}(48-7\lambda_{0})t^2+10\lambda_{0}(33-4\lambda_{0})t^2\nonumber\\
&+12(-3\lambda_{0}^2+36\lambda_{0}-8)t^2-6(5\lambda_{0}^2-54\lambda_{0}+24)t^2\nonumber\\
&-2(-25\lambda_{0}^2-36\lambda_{0}+36)t^2-24(1-\lambda_{0})\nonumber\\
=&(1806\lambda_{0}-312)t^2-24(1-\lambda_{0})>1830\lambda_{0}-336>0.
\end{align}
for $t\in(1,\sqrt[6]{2})$.

The inequality (2.18) implies that $G_{\lambda_{0}}(t)$ is strictly increasing in $[1,\sqrt[6]{2}]$. Then from (2.13) and (2.17) we clearly see that there exists $t_{1}\in(1,\sqrt[6]{2})$ such that $F_{\lambda_{0}}(t)$ is strictly decreasing in $[1,t_{1}]$ and strictly increasing in $[t_{1},\sqrt[6]{2}]$.

Therefore, $F_{\lambda_{0}}(t)<0$ for all $t\in(1,\sqrt[6]{2})$ follows from (2.11) and (2.16) together with the piecewise monotonicity of $F_{\lambda_{0}}(t)$. $\Box$

\medskip
\bigskip
\centerline{\bf 3. Proof of Theorems 1.1 and 1.2}
\setcounter{section}{3}\setcounter{equation}{0}
\bigskip

{\bf\em Proof of Theorem 1.1.} From (1.1) we clearly see that $M(a,b)$, $Q(a,b)$ and $A(a,b)$ are symmetric and homogeneous of degree $1$. Without loss of
generality, we assume that $a>b$. Let $p\in(0,1)$ and $x=(a-b)/(a+b)$, $t=\sqrt[6]{x^2+1}$ and $\alpha_{0}=(3-3\sqrt[6]{2}\log(1+\sqrt{2}))/[(2+\sqrt{2}-3\sqrt[6]{2})\log(1+\sqrt{2})]=0.777\cdots$. Then $x\in(0,1)$, $t\in(1,\sqrt[6]{2})$,
\begin{align}
&\frac{M(a,b)-Q^{1/3}(a,b)A^{2/3}(a,b)}{1/3Q(a,b)+2A(a,b)/3-Q^{1/3}(a,b)A^{2/3}(a,b)}\nonumber\\
=&\frac{3[x-\sqrt[6]{1+x^2}{\arcsinh}(x)]}{(\sqrt{1+x^2}-3\sqrt[6]{1+x^2}+2){\arcsinh}(x)}
\end{align}
and
\begin{align}
&p \left[\frac{1}{3}Q(a,b)+\frac{2}{3}A(a,b)\right]+(1-p)Q^{1/3}(a,b)A^{2/3}(a,b)-M(a,b)\nonumber\\
=&A(a,b)\left[p(\frac{1}{3}\sqrt{1+x^2}+\frac{2}{3})+(1-p)\sqrt[6]{1+x^2}-\frac{x}{{\arcsinh}(x)}\right]\nonumber\\
=&\frac{A(a,b)\left[p(\sqrt{1+x^2}+2)+3(1-p)\sqrt[6]{1+x^2}\right]}{3{\arcsinh}(x)}f_{p}(t).
\end{align}
where $f_{p}(t)$ is defined as in Lemma 2.1.
Note that
\begin{equation}
\lim\limits_{x\rightarrow
0}\frac{3[x-\sqrt[6]{1+x^2}{\arcsinh}(x)]}{(\sqrt{1+x^2}-3\sqrt[6]{1+x^2}+2){\arcsinh}(x)}=\frac{4}{5},
\end{equation}
\begin{equation}
\lim\limits_{x\rightarrow
1}\frac{3[x-\sqrt[6]{1+x^2}{\arcsinh}(x)]}{(\sqrt{1+x^2}-3\sqrt[6]{1+x^2}+2){\arcsinh}(x)}=\alpha_{0}.
\end{equation}

Therefore, Theorem 1.1 follows easily from (3.2)-(3.4) and Lemma 2.1. $\Box$

\bigskip
{\bf\em Proof of Theorem 1.2.} Since $M(a,b)$, $C(a,b)$ and $A(a,b)$ are symmetric and homogeneous of degree $1$. Without loss of
generality, we assume that $a>b$. Let $p\in(0,1)$ and $x=(a-b)/(a+b)$, $t=\sqrt[6]{x^2+1}$ and $\lambda_{0}=6(1-\sqrt[6]{2}\log(1+\sqrt{2}))/(7-6\sqrt[6]{2}\log(1+\sqrt{2}))=0.274\cdots$. Then $x\in(0,1)$, $t\in(1,\sqrt[6]{2})$,
\begin{align}
&\frac{M(a,b)-C^{1/6}(a,b)A^{5/6}(a,b)}{1/6C(a,b)+5A(a,b)/6-C^{1/6}(a,b)A^{5/6}(a,b)}\nonumber\\
=&\frac{6[x-\sqrt[6]{1+x^2}{\arcsinh}(x)]}{(x^2+6-6\sqrt[6]{1+x^2}){\arcsinh}(x)}
\end{align}
and
\begin{align}
&p \left[\frac{1}{6}C(a,b)+\frac{5}{6}A(a,b)\right]+(1-p)C^{1/6}(a,b)A^{5/6}(a,b)-M(a,b)\nonumber\\
=&A(a,b)\left[p\left(1+\frac{1}{6}x^2\right)+(1-p)\sqrt[6]{1+x^2}-\frac{x}{{\arcsinh}(x)}\right]\nonumber\\
=&\frac{A(a,b)\left[p(6+x^2)+6(1-p)\sqrt[6]{1+x^2}\right]}{6{\arcsinh}(x)}F_{p}(t).
\end{align}
where $F_{p}(t)$ is defined as in Lemma 2.2.
Note that
\begin{equation}
\lim\limits_{x\rightarrow
0}\frac{6[x-\sqrt[6]{1+x^2}{\arcsinh}(x)]}{(x^2+6-6\sqrt[6]{1+x^2}){\arcsinh}(x)}=\frac{8}{25},
\end{equation}
\begin{equation}
\lim\limits_{x\rightarrow
1}\frac{6[x-\sqrt[6]{1+x^2}{\arcsinh}(x)]}{(x^2+6-6\sqrt[6]{1+x^2}){\arcsinh}(x)}=\lambda_{0}.
\end{equation}

Therefore, Theorem 1.2 follows easily from (3.5)-(3.8) and Lemma 2.2. $\Box$

\medskip
{\bf REMARK 3.1}. If we take $\alpha=0$ and $\beta=1$ in Theorem 1.1, then the double inequality (1.4) reduces to (1.2).

\medskip
{\bf REMARK 3.2}. If we take $\lambda=0$ and $\mu=1$ in Theorem 1.2, then the double inequality (1.5) reduces to (1.3).

\medskip
\def\refname
\bigskip
{\small \hfill {\em Yuming Chu}

\hfill {\em Department of Mathematics}

\hfill {\em Huzhou Teachers College}

\hfill {\em Huzhou 313000}

\hfill {\em China}

\hfill {\em e-mail:} chuyuming@hutc.zj.cn

\bigskip
{\small \hfill {\em Miaokun Wang}

\hfill {\em Department of Mathematics}

\hfill {\em Huzhou Teachers College}

\hfill {\em Huzhou 313000}

\hfill {\em China}

\hfill {\em e-mail:} wmk000@126.com
\end{document}